\documentclass{article}
\usepackage{pdfpages}
\usepackage{graphicx}
\usepackage{graphicx}
\usepackage{amsmath}
\usepackage{subfigure}
\usepackage{float}
\usepackage{amssymb} \usepackage{amsthm}

\def\R{\hbox{{\rm I}\kern-0.2em{\rm R}\kern0.2em}}%mathematical R for reals
   
   \def\T{\cal T}

\def\be{\begin{equation}} \def\ee{\end{equation}} 
  \def\({\left(} \def\){\right)} 
\def\[{\left[}
\def\]{\right]}
\def\bc{\begin{center}}
\def\ec{\end{center}}

\begin{document}
\title{ Reductions and exact solutions of a \\cubic Schr\"odinger Partial Differential Equation }

\author{P. Masemola and T. Phidane}
\date{}
\maketitle
\begin{abstract}
Lie symmetry analysis is an established method for generating symmetries of differential equations. We apply this method together the generalized fundamental theorem of double reduction. In particular, Noether symmetries and some associated conservation laws are constructed in our investigation to find exact solutions of higher order partial differential equations and complex partial differential equations.
\end{abstract}

\section{Introduction}
In this paper we use the method of double reduction to find the exact solutions of a Cubic nonlinear Schr{\"o}dinger Partial Differential equation. Experiments on compressional dispersive Alven(CDA) waves have been done before on \cite{13} and \cite{14} where the relationship between CDA waves and the pertubations were analysed. Further in \cite{14}, the amplitude of the waves in a magnetic electron-positron plasma was discussed. In both experiments it was concluded that the system of equations under investigation can be written as the following equation:
\begin{equation}
u_{tt}-(3a^{2}+c^{2})u_{xx}-\delta^{2}u_{xxxx}-\delta^{2}u_{xxtt}=0.
\end{equation} 
This equation was analysed in \cite{12}.\newline
\\

In order to get an envelope of CDA waves we need an interaction of a CDA pump and a quasi stationary compressional magnetic field.When the CDA envelope evolves we get a cubic nonliniar Schr{\"o}dinger Partial Differential equation (NLSE) written as:
\begin{equation}\label{chr1}
iq_{t}+\beta iq_{x}-\gamma q_{xx}+\delta q |q|^{2}=0,
\end{equation}
where $q$ is the complex valued dependent variable, and $\gamma$ is the coefficient of the GVD, $\gamma$ is the self-phase modulation (SPM) owing to /kerr law $\beta$ is the inter-modal dispersion(IMD). We find that if we balance SPM and GVD, we get solitons.

To create a system of equation, we substitute $q = u+iv$ into \eqref{chr1} where $u=u(x,t)$, $v=v(x,t)$ and separate into real and imaginary parts to obtain:
\begin{align}
&G^{1}=u_{t}+\beta u_{x}-\gamma v_{xx}+\delta v(u^{2}+v^{2})=0,\\\nonumber
&G^{2}=-v_{t}-\beta v_{x}-\gamma u_{xx}+\delta(u^{2}+v^{2})=0.\label{chr2}
\end{align}

\newpage

\section{Conservation laws and conserved quantities of Cubic Schr{\"o}dinger Partial Differential equation}

To find multipliers and conservation laws, the condition below must hold.

$q^{1}(u_{t}+\beta u_{x}-\gamma v_{xx}+\delta v(u^{2}+v^{2}))+q^{2}(-v_{t}-\beta v_{x}-\gamma u_{xx}+\delta(u^{2}+v^{2})) = D_{t}T^{t}+D_{x}T^{x}$.\newline where $q^{1},q^{2}$ are multipliers of \eqref{chr1} and $T^{t},T^{x}$ conserved vectors of \eqref{chr1}.

From total divergence we set

$\frac{\delta}{\delta u}[q^{1}(u_{t}+\beta u_{x}-\gamma v_{xx}+\delta v(u^{2}+v^{2}))+q^{2}(-v_{t}-\beta v_{x}-\gamma u_{xx}+\delta(u^{2}+v^{2}))]= 0$ then use Maple and Mathematica to get the following results,

(i) ($q^{1},q^{2}$) = ($v_{x},u_{x}$)

\begin{align}
T^{x}_{1} =& \frac{1}{4}(\delta u^{4}+2\delta u^{2}v^{2}+\delta v^{4}+2vu_{t}-2uv_{t}-2\gamma (u_{x}^{2}+v_{x}^{2})),\\\nonumber
T^{t}_{1} =&\frac{1}{2}(-vu_{x}+uv_{x}).
\end{align}

(ii) ($q^{1},q^{2}$) = ($u,-v$)

\begin{align}\label{ii}
T^{x}_{2} =& \frac{1}{2}(\beta u^{2}+v(\beta v+2\gamma u_{x})-2\gamma uv_{x}),\nonumber\\
T^{t}_{2} =&\frac{1}{2}(u^{2}+v^{2}).
\end{align}

(iii) ($q^{1},q^{2}$) = ($u_{t},v_{t}$)

\begin{align}
T^{x}_{3} =& \frac{1}{2}(-\gamma(u_{t}u_{x}+v_{x}v_{t})+u(\beta v_{t}+\gamma u_{xt})+v(-\beta u_{t}+\gamma v_{xt})),\nonumber\\
T^{t}_{3} =&\frac{1}{4}(\delta u^{4}+2\delta u^{2}v^{2}-2u(\beta v_{x}+\gamma u_{xx})+v(\delta v^{3}+2\beta u_{x}-2\gamma v_{xx})).
\end{align}

(iv) ($q^{1},q^{2}$) = ($-\beta u+xu+2\gamma tv_{x},\beta t v-xv+2\gamma t u_{x}$)

\begin{align}
T^{x}_{4} =& \frac{1}{2}[t\gamma\delta u^{4}+\beta(x-t\beta)v^{2}+t\gamma\delta v^{4}+u^{2}(\beta(x-t\beta)+2t\gamma\delta v^{2})\nonumber\\
&2\gamma v(t u_{t}+(x-t\beta)u_{x})-2\gamma u(tv_{t}+(x-t\beta)v_{x})-2t\gamma^{2}(u_{x}^{2}+v_{x}^{2})],\nonumber\\
T^{t}_{4} =&\frac{1}{2}((x-t\beta)u^{2}+v((x-t\beta)v-2t\gamma u_{x})+2t\gamma uv_{x}).
\end{align}

\newpage

\section{Symmetries of Cubic Schr{\"o}dinger Partial Differential equation}
A one-parameter Lie group of transformations that leave \eqref{chr2} invariant will be written as a vector field
\begin{equation}
X=\tau(x,t,u,v)\partial_{t}+\xi(x,t,u,v)\partial_{x}+
\eta^{1}(x,t,u,v)\partial_{u}+\eta^{2}(x,t,u,v)\partial_{v}
\end{equation}
Equation \eqref{chr2} has the following variational symmetries:

\begin{align*}
X_{1}=&\partial_{t},\\
X_{2}=&\partial_{x},\\
X_{3}=&u\partial_{v}-v\partial_{u},\\
X_{4}=&2\gamma t\partial_{x}+(\beta t-x)u\partial_{v}+(-\beta t+x)v\partial_{u},\\
X_{5}=&2t\partial_{t}+(\beta t+x)\partial_{x}-u\partial_{u}-v\partial_{v}.
\end{align*}

\newpage

\section{Double reduction of Cubic Schr{\"o}dinger Partial Differential equation}

In order to perform a double reduction, we need to show association between the symmetries and the conservtion laws calculated in the previous sections.\newline
\\

To check if a Lie-B{\"a}cklund symmetry $X$ of differential equation is associated with conservation law$T=(T^{1},T^{2},...,T^{n})$ the following equation must hold
\begin{equation}\label{assoc}
T^{*i}=X(T^{i})+T^{i}D_{j}\xi^{j}-T^{j}D_{j}\xi^{i}=0
\end{equation} 
for $i = 1,2,...,n.$

\subsection*{Double reduction of Cubic Schr{\"o}dinger Partial Differential equation by symmeties $X_{1}$ and $X_{3}$}$ $\newline
In this section we use symmetries $<X_{1}$, $X_{3}>$ and conserved vectors given in \eqref{ii} to reduce the Cubic nonlinear Schr{\"o}dinger Partial Differential equation.\newline
\\
Firstly we show that Lie-B{\"a}cklund symmetries $X_{1}$ and $X_{3}$ are associated with conservation law $T=(T^{x},T^{t})$ given in \eqref{ii} using the following version of \eqref{assoc} for $i = 1,2:$

\begin{equation*}
T^{*}=X\begin{pmatrix}
T^{t}\\
T^{x}
\end{pmatrix} -\begin{pmatrix}
&D_{t}\xi^{t}&D_{x}\xi^{t}\\
&D_{t}\xi^{x}&D_{x}\xi^{x}
\end{pmatrix}\begin{pmatrix}
T^{t}\\T^{x}
\end{pmatrix}+(D_{t}\xi^{t}+D_{x}\xi^{x})\begin{pmatrix}
T^{t}\\T^{x}
\end{pmatrix}
\end{equation*}

We obtain\newline
\\

\begin{align*}
\begin{pmatrix}
T_{2}^{*t}\\
T_{2}^{*x}
\end{pmatrix}
 =& X_{1}^{[1]} \begin{pmatrix}
T_{2}^{t}\\T_{2}^{x}
\end{pmatrix} - \begin{pmatrix}
0&0\\
0&0
\end{pmatrix}
\begin{pmatrix}
T_{2}^{t}\\T_{2}^{x}
\end{pmatrix}+(0)\begin{pmatrix}
T_{2}^{t}\\T_{2}^{x}
\end{pmatrix}\\
=&X_{1}^{[1]}\begin{pmatrix}
&\frac{1}{2}(u^{2}+v^{2})\\
&\frac{1}{2}(\beta u^{2}+v(\beta v+2\gamma u_{x})-2\gamma uv_{x})
\end{pmatrix}\\
=&\begin{pmatrix}
0\\0
\end{pmatrix}
\end{align*}

and

\begin{align*}
\begin{pmatrix}
T_{2}^{*t}\\
T_{2}^{*x}
\end{pmatrix}
 =& X_{3}^{[1]} \begin{pmatrix}
T_{2}^{t}\\T_{2}^{x}
\end{pmatrix} - \begin{pmatrix}
0&0\\
0&0
\end{pmatrix}
\begin{pmatrix}
T_{2}^{t}\\T_{2}^{x}
\end{pmatrix}+(0)\begin{pmatrix}
T_{2}^{t}\\T_{2}^{x}
\end{pmatrix}\\
=&X_{3}^{[1]}\begin{pmatrix}
&\frac{1}{2}(u^{2}+v^{2})\\
&\frac{1}{2}(\beta u^{2}+v(\beta v+2\gamma u_{x})-2\gamma uv_{x})
\end{pmatrix}\\
=&\begin{pmatrix}
u\frac{1}{2}2v-v\frac{1}{2}2u\\
u\frac{1}{2}(2\beta v+2\gamma u_{x})-v\frac{1}{2}(2\beta u+2\gamma v_{x})+u_{x}\frac{1}{2}(-2\gamma u)-v_{x}(\frac{1}{2}2\gamma v)
\end{pmatrix}\\
=&\begin{pmatrix}
0\\0
\end{pmatrix}
\end{align*}
where
\begin{align*}
X_{1}^{[1]}=&\partial_{t}\\
X_{3}^{[1]}=&u\partial_{v}-v\partial_{u}+u_{t}\partial_{v_{t}}+u_{x}\partial_{v_{x}}-v_{t}\partial_{u_{t}}-v_{x}\partial_{u_{x}}
\end{align*}
Thus $X_{1}$ and $X_{3}$ are both associated with $T_{2}$. We then consider a linear combination of $X_{1}$ and $X_{3}$, of the form $X=X_{1}+cX_{3}$ and transform this generator to its canonical form $Y$ = $\frac{\partial}{\partial_{s}}$, where we assume that this generator is of the form

\begin{equation*}
Y=0\frac{\partial}{\partial_{r}}+\frac{\partial}{\partial_{s}}+0\frac{\partial}{\partial_{w}}+0\frac{\partial}{\partial_{p}}.
\end{equation*}

From $X(r)=0$, $X(s)=1$, $X(w)=0$, and $X(p)=0$, we obtain

\begin{equation}\label{nne}
\frac{dt}{1}=\frac{dx}{0}=\frac{du}{-cv}=\frac{ds}{1}=\frac{dr}{0}=\frac{dw}{0}=\frac{dp}{0}=\frac{dv}{cu}.
\end{equation}

We then solve \eqref{nne} using the method of invariants, and the results are summarized in Table 4.1 

\begin{table}[h!]
\centering
\begin{tabular}{||c  c||} 
 \hline
 &\\
\multicolumn{2}{||c||}{Invariats of $X=\frac{\partial}{\partial_{t}}+c(u\partial_{v}-v\partial_{u})$}\\

&\\
\hline\hline
 &\\
 $\frac{dt}{1}=\frac{ds}{1}$ & $b_{1}$ = $s-t $\\ 
 &\\
 $\frac{du}{-cv}=\frac{dv}{ku}$ & $b_{2}=u^{2}+v^{2}$\\
 &\\
 $\frac{dv}{cu}$ = $\frac{dt}{1}$ & $b_{3} = \arctan(\frac{v}{u})-ct$\\
 &\\
 $\frac{dr}{0}$&$b_{4}=r$\\
 &\\
 $\frac{dp}{0}$&$b_{5}=p$\\
 &\\
$\frac{dw}{0}$ & $b_{6}$ = $w$\\
&\\
$\frac{dx}{0}$ &$b_{7}$ =$x$ \\ [1ex] 
 \hline
\end{tabular}
\caption{invariants table}
\label{table:1}
\end{table}
.

\begin{align*}
&&
\end{align*}

\newpage

\begin{align*}
&&
\end{align*}
By choosing $b_{1}=0$, $b_{4}=b_{7}$, $b_{6} =\sqrt{b_{2}}$,and $b_{3}=b_{5}$ we get
\begin{align*}
&r = x\\
&s = t\\
&w = \sqrt{u^{2}+v^{2}}\\
&p = \arctan\big(\frac{v}{u}\big)-ct
\end{align*} 

The matrices $A$ and $A^{-1}$ can be computed using the canonical coordinates above

\begin{equation}
A = \begin{pmatrix} 
D_{s}t & D_{s}x \\
D_{r}t & D_{r}x \\
\end{pmatrix}=\begin{pmatrix} 
1 & 0 \\
0 & 1 \\
\end{pmatrix}=(A^{-1})^{T}
\end{equation}
and $J$ = $\det(A)$ = 1.

The inverse canonical coordinates are presented below
\begin{align*}
&x = r,\\
&t = s,\\
&u = w\cos(p+cs),\\
&v = w\sin(p+cs).\\
\end{align*}
The first and second partial derivatives of $u$ and $v$ in terms of new dependent variables $w$ and $p$ are,

\begin{align*}
u_{x} =&w_{r}\cos(p+cs)-wp_{r}\sin(p+cs),\\
u_{xx} =&w_{rr}\cos(p+cs)-2w_{r}p_{r}\sin(p+cs) -wp_{r}^{2}\cos(p+cs)-wp_{rr}\sin(p+cs),\\
u_{t} =&-cw\sin(p+cs),\\
v_{x} =&w_{r}\sin(p+cs)+wp_{r}\cos(p+cs),\\
v_{xx} =&w_{rr}\sin(p+cs)+2w_{r}p_{r}\cos(p+cs) +wp_{rr}\cos(p+cs)-wp_{r}^{2}\sin(p+cs),\\
v_{t} =&cw\cos(p+cs).
\end{align*}
The reduced conserved form will be

\begin{align}\label{ndibu}
\begin{pmatrix}
T_{2}^{s}\\
&\\
T_{2}^{r}\\
\end{pmatrix}
=&J(A^{-1})^{T}\begin{pmatrix}
T_{2}^{t}\\
&
T_{2}^{x}\\
\end{pmatrix}\nonumber\\
=&\begin{pmatrix}
&\frac{1}{2}(u^{2}+v^{2})\\
&\\
&\frac{1}{2}(\beta u^{2}+v(\beta v+2\gamma u_{x})-2\gamma uv_{x})
\end{pmatrix}.
\end{align}

Substituting the first and second partial derivatives of $u$ and $v$ into \eqref{ndibu} we get
\begin{align*}\label{ndibu1}
\begin{pmatrix}
T_{2}^{s}\\
&\\
T_{2}^{r}\\
\end{pmatrix}=\begin{pmatrix}
&\frac{1}{2}w(r)^{2}\\
&\\
&\frac{1}{2}(2\beta w(r)^{2}+\gamma w(r)w(r)_{r}\sin(2(p+ks))+2\gamma w(r)^{2}\cos(2(p+ks))
\end{pmatrix}
\end{align*}

where the reduced conserved form is also given by
\begin{equation}
D_{s}T_{2}^{s}=0.
\end{equation}

The second step of double reduction can also be given as
\begin{equation}
w(r)^{2}=\epsilon
\end{equation}
$\epsilon$ is constant, or equivalently

\begin{equation}\label{poa}
w(r)=\sqrt{\epsilon}
\end{equation}
Differentiating \eqref{poa} implicitly with respect to $r$ we get
\begin{equation}\label{poa1}
\frac{d}{dr}w(r)=0.
\end{equation}

Our original system is equivalent to
\begin{equation}\label{d5}
sys_{1}=
\left\{%
\begin{array}{l}
    \hbox{${g_1}^{1}Q^{1}+{ g_1}^{2}Q^{2}=0,$} \\
    \hbox{${ g_1}^{1}Q^{1}-{ g_1}^{2}Q^{2}=0.$} \\
\end{array}%
\right.
\end{equation}

This system can be rewritten as
\begin{eqnarray}\label{d6}
D_t {T_{2}}^t+D_x { T_{2}}^x &=&0, \nonumber\\[1\jot]
{ G_{1}}^{1}Q^{1}-{ G_{1}}^{2}Q^{2} &=&0.
\end{eqnarray}
The second equation of $sys_{1}$ from is given by

\begin{equation}\label{poa2}
u(u_{t}+\beta u_{x}-\gamma v_{xx}+\delta v(u^{2}+v^{2}))-(-v)(-v_{t}-\beta v_{x}-\gamma u_{xx}+\delta(u^{2}+v^{2}) = 0
\end{equation}
Substituting \eqref{poa},\eqref{poa1},\eqref{poa2} and the first and second partial derivatives of $u$ and $v$ in terms of $w(r)$ and $p(r)$ we get

\begin{align}\label{po3}
&-c\epsilon\sin(2p(r)+2cs)-\beta\epsilon p'(r)\sin(2p(r)+2cs)-\gamma\epsilon p''(r)\cos(2p(r)+2cs)\nonumber\\&+\gamma\epsilon p'(r)^{2}\sin(2p(r)+2cs)+\delta\epsilon^{2}\sin(2p(r)+2cs)=0
\end{align}
 
When computing the final solution to equation \eqref{po3} we get tedious solution. We then compute the solution of \eqref{po3} for the following cases.\newline
\\

Case 1: $c=\gamma=0$

\begin{align}
p(r)=&0\label{weqa11},\\
p(r)=&-\frac{\pi}{2}\label{weqa22},\\
p(r)=&\frac{\pi}{2}\label{weqa33},\\
p(r)=&\frac{r\delta\epsilon}{\beta}+c_{1}.\label{weqa44}
\end{align}

From \eqref{weqa11}, $p(r)=0$:

\begin{align*}
u=\sqrt{\epsilon},\\
v=0,\\
q=\sqrt{\epsilon}.
\end{align*}

From \eqref{weqa22}, $p(r)=-\frac{\pi}{2}$

\begin{align*}
u=0,\\
v=-\sqrt{\epsilon},\\
q=-i\sqrt{\epsilon}.
\end{align*}

From \eqref{weqa33}, $p(r)=\frac{\pi}{2}$

\begin{align*}
u=0,\\
v=\sqrt{\epsilon},\\
q=i\sqrt{\epsilon}.
\end{align*}

From \eqref{weqa44}, $p(r)=\frac{r\delta\epsilon}{\beta}+c_{1}$

\begin{align*}
u=\sqrt{\epsilon}\cos(\frac{x\delta\epsilon}{\beta}+c_{1}),\\
v=\sqrt{\epsilon}\sin(\frac{x\delta\epsilon}{\beta}+c_{1})
.
\end{align*}
From canonical coordinates
\begin{equation}
q=\sqrt{\epsilon}e^{i(\frac{x\delta\epsilon}{\beta}+c_{1})}
\end{equation}

Case 2 $c=\beta=0$:

\begin{align}
p(r)=&0,\\
p(r)=&-\frac{\pi}{2},\\
p(r)=&\frac{\pi}{2},\\
p(r)=&c_{1}.\label{wasaq}
\end{align}

From \eqref{wasaq} $p(r)=c_{1}$

\begin{align*}
u=\sqrt{\epsilon}\cos(c_{1}),\\
v=\sqrt{\epsilon}\sin(c_{1}),\\
q=\sqrt{\epsilon}e^{ic_{1}}.
\end{align*}

Case 3 $\delta=\gamma=0$:

\begin{align}
p(r)=&-cs\label{poi1},\\
p(r)=&-\frac{\pi}{2}-cs,\\
p(r)=&\frac{\pi}{2}-cs,\\
p(r)=&-\frac{cr}{\beta}+c_{1}.
\end{align}

From \eqref{poi1} $p(r)=-\frac{cr}{\beta}+c_{1}$

\begin{align*}
u=\sqrt{\epsilon}\cos(\frac{c(\beta s- r)}{\beta}+c_{1})\\
v=\sqrt{\epsilon}\sin(\frac{c(\beta s- r)}{\beta}+c_{1})\\
\end{align*}
From canonical coordinates
\begin{equation}
q=\sqrt{\epsilon}e^{i(\frac{c(\beta t- x)}{\beta}+c_{1})}
\end{equation}

\newpage

\section{Conclusion}
In this paper we have calculated multipliers and conservation laws of the respective Cubic Schr\"odinger partial differential equation. Further, we showed that if a symmetry is associated with a conservation law we can apply the fundamental theorem of double reduction obtain an ordinary differential equation. Consequently we used Maple and Mathematica to calculate exact solutions for special cases.

\section{Acknowledements}
This project is supported by the NRF Thuthuka under grant number TTK14051567369 and without this, the project would not be possible.

\end{document}